\documentclass{article}
\usepackage{graphicx}
\usepackage{amsmath}
\usepackage{amsfonts}
\usepackage{amssymb}
\newtheorem{theorem}{Theorem}

\newtheorem{example}[theorem]{Example}

\newtheorem{proposition}[theorem]{Proposition}

\begin{document}
\title{A constructive version of the Boyle-Handelman
theorem on the spectra of nonnegative
matrices}
\maketitle


\begin{center}
 Thomas
J. Laffey\ \ 

 School
of Mathematical Sciences

 University
College Dublin

Belfield,

Dublin 4, Ireland

E-mail:
Thomas.Laffey@ucd.ie\textbf{\ \ \ }
\end{center}

\bigskip\medskip

\bigskip

\begin{abstract}
 
\smallskip

A constructive version of the celebrated Boyle-Handelman theorem on the non-zero
spectra of nonnegative matrices is presented.


\bigskip
\end{abstract}

\section{Introduction.}

\ Let%
\[
\sigma=(\lambda_{1},\quad...\quad,\lambda_{n})
\]
be a list of complex numbers \ and let%
\[
s_{k}:=\lambda_{1}^{k}+\quad...\quad+\lambda_{n}^{k},\text{ }k=1,2,3,...
\]
The nonnegative inverse eigenvalue problem (\textbf{NIEP}) asks for necessary
and sufficient conditions on $\sigma$ in order that it be the spectrum of an entry-wise nonnegative
matrix. If this occurs, we say that $\sigma$ is \textit{realizable, }and we call a nonnegative matrix $A$
with spectrum $\sigma$ a \textit{realizing matrix} for $\sigma$.

A necessary condition for realizability coming from the Perron-Frobenius
theorem [2] is that there exists $j$ with $\lambda_{j}$ real and $\lambda_{j}\geq\mid\lambda
_{i}\mid$, for all $i.$ \ Such a $\lambda_{j}$ is called the \textit{Perron
root} of $\sigma$. 

A more obvious necessary condition is that all the $s_{k}$ are nonnegative.
A stronger form of this condition was found independently by Loewy and London [11]
and Johnson [8], namely:%
\[
\text{(\textbf{JLL}) }n^{k-1}s_{km}\geq s_{m}^{k}\text{, for all positive
integers, }k\text{ and }m.
\]
In terms of $n$, a complete solution of the NIEP is available only for
$n\leq4$. \ The solution for $n=4$, expressed in terms of inequalities for the $s_{k}$, appears in the
\ PhD thesis of Meehan [12] and a solution in terms of the coefficients of the characteristic
polynomial has been published more recently by Torre-Mayo, Abril-Raymundo, Alarcia-Estevez, Marijuan, and Pisanero [14].

However, the same problem in which we may augment the list $\sigma$ by adding
an arbitrary number $N$ of zeros was solved by Boyle and Handelman [4]. Using a range of tools coming from linear algebra, dynamical systems,ergodic theory, and graph theory, they proved the remarkable result that if
\begin{enumerate}
\item $\sigma$ has  a Perron element $\lambda_{1}>\mid\lambda_{j}\mid$(all $j>1$) \ 
and
\item $s_{k}\geq0$ \ for all positive integers $k$ \ 
(and $s_{m}=0$ for some $m$ implies $s_{d}=0$ \ for all positive divisors $d$
of $m$),
then
\[
\sigma_{N}:=(\lambda_{1},\quad...\quad,\lambda_{n},0.\quad...\quad,0)\text{
}(N\text{ zeros})
\]
is realizable for all sufficiently large $N$.
\end{enumerate}
Under these assumptions, a realizing matrix can be chosen to be primitive.
\ See Friedland [6] for an extension to the irreducible case. \ \ \ \ \ \ \ \ \ \ \ 

The proof of the Boyle-Handelman result is not constructive and does not
provide a bound on the minimal number $N=N(\sigma)$of zeros required for realizability.

Finding a constructive proof, with a bound on the minimum number$N$ of zeros
required, has been an area of much research, and a number of special cases have been resolved. In
particular, a best possible result in the case that Re($\lambda_{j}$) $\leq0$, for all $\ j>1, $
has been obtained by \v{S}migoc and the author [9] and, when $\sigma$ is real and has exactly two positive
entries, a constructive proof with a bound on $N$ has also been found [10].

In the case that $\sigma$ is real and has just one positive entry, then the
inequality $s_{1}\geq0$ is necessary and sufficient for realizability. This was proved by Suleimanova [13] and this
is often viewed as the first result on the NIEP. \ Friedland [5] re-proved her result by showing
that the companion matrix with spectrum $\sigma$ has nonnegative entries, and matrices
related to companion matrices are used in the cited work with \v{S}migoc.

Here, a constructive approach to the Boyle-Handelman result is presented. It
is shown that a certain kind of patterned matrix is ''universal'' for the realization of
spectra with power sums$s_{k}>0$, ( $k=1,2,3,$ \ \ $...$ \ \ ). in the sense that all such spectra
satisfying the Perron condition (i) above can, with sufficiently many zeros added, be realized as
the spectrum of a primitive nonnegative matrix with that pattern.

\medskip

\bigskip

\section{A matrix related to Newton's identities}

\smallskip

\bigskip

Let
\begin{align*}
\tau &  =(\mu_{1},\quad...\quad,\mu_{n}),\\
x_{k}  &  :=\mu_{1}^{k}+\quad...\quad+\mu_{n}^{k},\text{ }k=1,2,3,...\\
q(x)  &  :=\Pi_{i=1}^{n}(x-\mu_{i})\\
&  =x^{n}+q_{1}x^{n-1}+\quad...\quad+q_{n}.
\end{align*}

\bigskip

Let $X_{n}=$%

\[
\left(
\begin{array}
[c]{cccccccccc}%
x_{1} & 1 & 0 &  & . & . & . &  &  & 0\\
x_{2} & x_{1} & 2 & 0 &  & . & . & . &  & 0\\
x_{3} & x_{2} & x_{1} & 3 & 0 &  & . & . & . & 0\\
. & x_{3} & . & . & . & . &  &  &  & \\
. & . &  & . & . & . &  &  &  & \\
. &  &  &  & . &  &  &  &  & \\
&  &  &  &  & . &  &  &  & \\
&  &  &  &  &  &  &  &  & \\
x_{n-1} & x_{n-2} &  & . & . & . &  & x_{2} & x_{1} & n-1\\
x_{n} & x_{n-1} &  &  & . & . & . & x_{3} & x_{2} & x_{1}%
\end{array}
\right)
\]

\bigskip

The matrix $X_{n}$ occurs in the context of the Newton identities relating the
coefficients of a polynomial to the power sums of its roots. If we use Cramer's rule to
express the $q_{i}$ in terms of the $x_{j}$, we get
\[
\det(X_{n})=(-1)^{n}n!q_{n}.
\]
However, the matrix $X_{n}$ itself, as distinct from its determinant, does not
appear to have been widely investigated. A key observation is:

\medskip

\begin{proposition}
 The characteristic polynomial of $X_{n}$ is%
\[
Q(x)=x^{n}+nq_{1}x^{n-1}+n(n-1)q_{2}x^{n-2}+\quad...\quad+n!q_{n}.
\]
\medskip
Since $X_{n}$ has nonnegative entries if the $x_{i}$ are nonnegative, it
follows that the spectrum of $Q(x)$ is realizable if the $x_{i}$, $(i=1,2,\quad...\hspace{0in}\quad,n),$ are nonnegative.
\end{proposition}
\hspace*{0in}

Suppose that we are given a list $\sigma=(\lambda_{1},\quad...\quad,\lambda_{n})$ that we wish to realize as the spectrum of a nonnegative matrix.

Let%
\begin{align*}  
f(x)  &  :=\Pi_{i=1}^{n}(x-\lambda_{i})\\
&  =x^{n}+p_{1}x^{n-1}+\quad...\quad+p_{n}.
\end{align*}
Let%
\[
q(x):=x^{n}+q_{1}x^{n-1}+\quad...\quad+q_{n}%
\]
where%
\[
q_{i}=\frac{p_{i}}{n(n-1)\quad...\quad(n-i+1)}\text{, }i=1,2,\quad
...\quad,n\text{.}%
\]

Then the corresponding $Q(x)$ is $f(x).$ Now the power sums $x_{i}$ of the
roots of $q(x)$ are nonnegative if and only if that holds for%
\[
x^{n}+nq_{1}x^{n-1}+\quad...\quad+n^{n}q_{n}.
\]
Hence we have

\medskip

\begin{theorem} \ $\sigma$ is realizable by the matrix $X_{n}$ if the $j $
th power sum of the roots of the polynomial%
\[
J_{n}(f(x)):=x^{n}+p_{1}x^{n-1}+\frac{n}{n-1}p_{2}x^{n-2}+\frac{n^{2}%
}{(n-1)(n-2)}p_{3}x^{n-3}+\quad...+\frac{n^{n-1}}{(n-1)!}p_{n}%
\]
is nonnegative for $j=1,2,3,\quad...\quad,n.$
\end{theorem}

\medskip

But now suppose that we choose $N>n$ and ask for the realizability of $\sigma$
with $N-n$ zeros added. This amounts to replacing $f(x)$ by $x^{N-n}f(x)$ and $J_{n}(f(x))$ by
$x^{N-n}J_{N}(f(x))$, where%

\begin{align}
J_{N}(f(x))& :=x^{n}+p_{1}x^{n-1}+\frac{N}{N-1}p_{2}x^{n-2} \\
&+\frac{N^{2}
}{(N-1)(N-2)}p_{3}x^{n-3}+...+\frac{N^{n-1}}{(N-1)(N-2)...(N-n+1)}p_{n}.
\end{align}

So $\sigma$ with $N-n$ zeros added is realizable by the matrix $X_{N}$ if the
$j$ th power sums of the roots of the polynomial $J_{N}(f(x))$ are nonnegative for $j=1,2,3,\quad...\quad,N.$

But observe that as $N\rightarrow\infty,J_{N}(f(x))\rightarrow f(x)$, since
$n$ is fixed.

Suppose that the power sums $s_{j}$ of the elements of $\sigma$ are positive
for all $j\geq1.$ Then, on continuity grounds, one might expect that for sufficiently large $N$, the
power sums of the roots of $J_{N}(f(x))$ would also be positive. However, this is not true in
general, but it is  true if $\sigma$ has its Perron element

$\ \ \ \ \ \ \ \ \ \ \ \ \ \ \ \ \lambda_{1}>\mid\lambda_{j}\mid
(j=1,2,\quad...\quad,n).$

In this case, $\sigma$ with sufficiently many zeros added is the spectrum of a
nonnegative matrix $X_{N}.$

Since, we only require that the $j$ th power sum of the roots of $J_{N}(f(x))
$ be nonnegative for $j=1,2,$ \ \ ... $\ \ ,N$, one can obtain a bound on the minimal number of
zeros required.

\medskip

\hspace*{0in} \ \ 

\section{ Main Theorem}

We now state the main result of this paper.

\smallskip

\begin{theorem}
Let
\[
\sigma=(\lambda_{1},\quad...\quad,\lambda_{n}),
\]
be a list of complex numbers with corresponding power sums%
\[
s_{k}:=\lambda_{1}^{k}+\quad...\quad+\lambda_{n}^{k},\quad k=1,2,3,...\text{
\ \ .}%
\]
Suppose that

(i) \ $\lambda_{1}>\mid\lambda_{j}\mid$, \ (all $j>1$) \ \ \ 

(ii) \ \ $s_{1}\geq0$, \ and\ $s_{m}>0$, \ for all $m\geq2$.

Let
\begin{align*}
f(x)  &  =\Pi_{i=1}^{n}(x-\lambda_{i})\\
&  =x^{n}+p_{1}x^{n-1}+\quad...\quad+p_{n}.
\end{align*}%

\[
\gamma=2\max(1,\mid p_{1}\mid,\mid p_{2}\mid^{1/2},\quad...\quad,\mid
p_{n}\mid^{1/n}).
\]%

\begin{align*}
\lambda_{0}  &  =\max\{\mid\lambda_{j}\mid:j>1\},\\
R  &  =\frac{(\lambda_{1}-\lambda_{0})}{4},\quad\ell=\frac{3\lambda
_{1}+\lambda_{0}}{\lambda_{1}+3\lambda_{0}},\quad r=\min(R,1),
\end{align*}%

\[
m=\max\{1,\lambda_{1}\},\quad N_{0}=\left\lceil \frac{\ln(2n-2)}{\ln(\ell
)}\right\rceil ,
\]%

\[
M=\min\{1,s_{2},\quad...\quad,s_{N_{0}}),
\]

and
\[
N=\left\lceil 2\left(  \frac{16\gamma nN_{0}(m+r)_{{}}^{N_{0}-1})}{3^{1/2}%
Mr}\right)  ^{n}\right\rceil \text{.}%
\]

Then $\sigma$ with $N-n$ zeros added is the spectrum of the nonnegative matrix
$X_{N}$, with $x_{k}:=\mu_{1}^{k}+\quad...\quad+\mu_{n}^{k},k=1,2,3,\quad...$ \ $,N,$ \ where%

\[
J_{N}(f(x))=(x-N\mu_{1})(x-N\mu_{2})...(X-N\mu_{n}).
\]
\end{theorem}
\quad

\bigskip

Given a list $\sigma$ satisfying the hypotheses, it is relatively easy to find
$N$ for which $J_{N}(f(x))$ has the corresponding power sums nonnegative, so one obtains a reasonably efficient constructive algorithm. However, the number of zeros required in the construction is not
optimal in general.

\bigskip

\section{Proofs of the results}

\medskip

{\small Let} $P=$

$\left(
\begin{array}
[c]{cccccccccc}%
1 & 0 & 0 &  & . & . & . &  &  & 0\\
q_{1} & 1 & 0 & 0 &  & . & . & . &  & 0\\
q_{2} & q_{1} & \frac{1}{2} & 0 & 0 &  & . & . & . & 0\\
q_{3} & q_{2} & \frac{q_{1}}{2} & \frac{1}{6} & 0 & . &  &  &  & \\
. & . & \frac{q_{2}}{2} & \frac{q_{1}}{6} & \frac{1}{24} & . &  &  &  & \\
. &  &  & \frac{q_{2}}{6} & . & . &  &  &  & \\
&  &  &  &  & . & . &  &  & \\
q_{n-2} &  &  &  &  &  &  & . &  & \\
q_{n-1} & q_{n-2} & . & . & . & . &  & . & \frac{1}{(n-2)!} & 0\\
q_{n} & q_{n-1} & \frac{q_{n-2}}{2} &  & . & . & . & . & \frac{q_{1}}{(n-2)!}%
& \frac{1}{(n-1)!}%
\end{array}
\right)  $

\medskip

and let $C=$

$\left(
\begin{array}
[c]{cccccccccc}%
0 & 1 & 0 &  & . & . & . &  &  & 0\\
0 & 0 & 1 & 0 &  & . & . & . &  & 0\\
._{{}} & . & 0 & 1 & 0 &  & . & . & . & 0\\
. &  & . & . & . & . &  &  &  & \\
. & . &  & . & . & . &  &  &  & \\
. &  &  &  & . &  &  &  &  & \\
&  &  &  &  & . &  &  &  & \\
&  &  &  &  &  &  &  &  & \\
0 & 0 &  & . & . & . &  & 0 & 0 & 1\\
-n!q_{n} & . &  &  & . & . & . & . & . & -nq_{1}%
\end{array}
\right)  $

be the companion matrix of $Q(x)=x^{n}+nq_{1}x^{n-1}+n(n-1)q_{2}x^{n-2}%
+\quad...\quad+n!q_{n}.$

Direct multiplication, using the Newton identities, yields $PC=X_{n}P$. This
proves the proposition.

\ 

To obtain the desired bound we use the following refinement by Bhatia, Elsner
and Krause [3 ] of a classical result of Ostrowski.

\begin{theorem} Let $f(x)=x^{n}+a_{1}x^{n-1}+\ \ ...\ \ +a_{n}$ and

\ \ \ \ \ \ \ \ \ \ \ \ \ \ \ \ \ \ \ \ \ \ \ $g(x)=x^{n}+b_{1}x^{n-1}%
+\ \ ...\ \ +b_{n}$

be real polynomials with roots $\alpha_{1},$ \ \ $...$ \ \ $,\alpha_{n}$ and
$\beta_{1},$ \ \ $...$ \ \ $,\beta_{n}$, respectively.

Then there is a labelling of $\beta_{1},$ \ \ $...$ \ \ $,\beta_{n}$ such that%

\[
\max\{\mid\alpha_{i}-\beta_{i}\mid:1\leq i\leq n\}\leq\left(  \frac{16}%
{3\sqrt{3}}\right)  (\sum{}_{k=1}^{n}\mid a_{k}-b_{k}\mid\gamma^{n-k}%
)^{1/n}\text{,}%
\]
where $\gamma=2\max\{\mid a_{k}\mid^{1/k},\mid b_{k}\mid^{1/k}:1\leq k\leq n\}.$
\end{theorem}

[The original Ostrowski result had the factor $(2n-1)$ in place of $\left(
\frac{16}{3\sqrt{3}}\right)  $].

Now let
\begin{align*}
f(x)  &  =(x-\lambda_{1})\text{ \ \ \ }....\text{ \ \ }(x-\lambda_{n})\\
&  =x^{n}+p_{1}x^{n-1}+\text{ \ \ }...\text{ \ \ }+p_{n}%
\end{align*}
and%
\[
g(x)=x^{n}+p_{1}x^{n-1}+\left(  \frac{N}{N-1}\right)  p_{2}x^{n-2}%
+\ \ ...\ \ +\left(  \frac{N^{n-1}}{(N-1)...(N-n+1)}\right)  p_{n}.
\]

We note that if $g(x)$ has nonnegative Newton power sums, then the
corresponding matrix $X_{N}$ is nonnegative and has spectrum $N\lambda_{1},$ \ \ ... \ \ \ $,N\lambda_{n}.$

Suppose that $\lambda_{1}>\mid\lambda_{j}\mid$ \ (all $j>1$) and let
$\lambda_{0}=\max(\mid\lambda_{j}\mid:j=2,$ \ \ ... \ \ $,n\}$ and
$R=\frac{\lambda_{1}-\lambda_{0}}{4}.$ \ Let $\ell=\frac{\lambda_{0}%
+R}{\lambda_{1}-R}$, so $\ell<1$. Let $r=\min\{1,R\}.$
Let%
\[
s_{k}=\lambda_{1}^{k}+\text{ \ \ ... \ \ }+\lambda_{n}^{k}%
\]
for $k=1,2,$ \ \ $...$ \ \ . Assume that $s_{1}\geq0,$ and that $s_{k}>0$, for
all $k>1.$
Let
\[
M=\min\{1,s_{k}:k=2,3,...\}.
\]
Let $\mu_{1},$ \ \ ... \ \ \ $,\mu_{n}$ be the roots of $g(x)=0$ and suppose that%
\[
\max\{\mid\lambda_{j}-\mu_{j}\mid:j=1,2,\text{ \ \ ... \ \ },n\}<\delta,\text{
\ \ (*)}%
\]
where%
\[
\delta=\frac{Mr}{nN_{0}(m+r)^{N_{0}-1}},
\]
with
\[
m=\max\{1,\lambda_{1}\},\quad N_{0}=\left\lceil \frac{\ln(2(n-1))}{\ln
(1/\ell)}\right\rceil \text{.}%
\]
Then $\mid\mu_{1}\mid$ is greatest among all the $\mid\mu_{j}\mid$, and ,
since $g(x)$ has real coefficients, $\mu_{1}$ is real and, since $\lambda_{1}$ is positive, so is $\mu_{1}$.
Let%
\[
S_{k}=\mu_{1}^{k}+\text{ \ \ ... \ \ }+\mu_{n}^{k}.
\]
Then $\mid s_{k}-S_{k}\mid\leq\sum$ $_{i=1}^{n}\mid\lambda_{i}^{k}-\mu_{i}%
^{k}\mid.$ Now%
\begin{align*}
&  \mid\lambda_{i}^{k}-\mu_{i}^{k}\mid=\mid\lambda_{i}^{{}}-\mu_{i}^{{}}%
\mid\mid\lambda_{i}^{k-1}+\lambda_{i}^{k-2}\mu_{i}+\text{ \ \ ...
\ }+\text{\ }\mu_{i}^{k-1}\mid\\
&  <\delta k(\lambda_{1}+r)^{k-1}.
\end{align*}
Suppose that $k\geq N_{0}$.
Then $S_{k}\geq(\lambda_{1}-r)^{k}-(n-1)(\lambda_{0}+r)^{k}=(\lambda
_{1}-r)^{k}(1-(n-1)(\frac{(\lambda_{0}+r)}{(\lambda_{1}-r)})^{k})>(\frac{1}%
{2})(\lambda_{1}-r)^{k}>0.$
For $k\leq N_{0}$,
\begin{align*}
&  \mid s_{k}-S_{k}\mid\leq\delta(1+2(\lambda_{1}+r)+\text{ \ \ ...
\ \ }+N_{0}(\lambda_{1}+r)^{N_{0}-1})\\
&  <\delta nN_{0}(m+r{\Large )}^{N_{0}-1}=Mr\leq M.
\end{align*}

So $S_{k}\geq0$, for all $k\geq2.$ \ Also, $S_{1}=s_{1}\geq0$. \ This shows
that if we can choose $N$ so that the inequality%
\[
\max\{\mid\lambda_{j}-\mu_{j}\mid:j=1,2,\text{ \ \ ... \ \ },n\}<\delta
\]
holds for that $\delta$, then the corresponding $X_{N}$ will be a nonnegative
matrix with spectrum
$\lambda_{1},$ \ \ ... \ \ \ $,\lambda_{n}$ and $N-n$ zeros.
Now,
\begin{align*}
\max&\{\mid\lambda_{j}-\mu_{j}\mid:j=1,2,...,n\} \leq \\
&(\frac
{16}{3\sqrt{3}})(\sum{}_{k=1}^{n}\mid p_{k}(\frac{N^{k-1}}{(N-1)\text{
\ \ ...}^{{}}(N-k+1)}-1)\gamma^{n-k}\mid^{1/n}.
\end{align*}

But%
\[
\frac{N^{k-1}}{(N-1)\text{ \ \ ...}^{{}}(N-k+1)}-1\leq\frac{2n^{2}}%
{N},\text{if }N>n^{2}.
\]
By definition, $\gamma=2\max_{{}}\{1,\mid p_{k}\mid^{1/k},k=1,2,$ $\ \ ...$
$\ \ ,n\}$. \ Hence
\[
\max\{\mid\lambda_{j}-\mu_{j}\mid:j=1,2,\text{ \ \ ... \ \ },n\}\leq
(\frac{16\gamma}{3\sqrt{3}})(\frac{2n^{3}}{N})^{1/n}.
\]
But $n^{1/n}\leq3^{1/3}.$ \ Hence%
\[
\max\{\mid\lambda_{j}-\mu_{j}\mid:j=1,2,\text{ \ \ ... \ \ },n\}\leq
\frac{16\gamma.2^{1/n}.}{\sqrt{3}N^{1/n}}\leq\delta\text{,}%
\]
provided%
\[
N\geq\frac{2^{4n+1}\gamma^{n}}{3^{n/2}\delta^{n}}=\frac{2(16\gamma
nN_{0}(m+r)^{N_{0}-1})^{n}}{3^{n/2}M^{n}r^{n}}\text{.}%
\]
This gives the required bound.\bigskip

There are variations of the Ostrowski bound, some using the Bombieri norm in
place of the $\ell_{2}$ one, available though the work of Beauzamy [1], Galantai and Hegedus [7],
\ and these may lead to better bounds for $N$ in certain circumstances. However, the main
interest is that such a bound exists, and the general form it has.

When the Perron root $\lambda_{1}=1,$a nonnegative matrix $A$ with the given
nonzero spectrum can be made stochastic. In this case $r$ and $\ell$ are measures of the spectral gap,
which control the rate at which the powers of $A$ converge to the stationary state of the
corresponding Markov process. The size of $N_{0}$ is inversely related to $r$ and $\ell.$

The number $M$ measures how close to zero the power sums can get, and we see
its appearance (as $M^{n}$) in the denominator of the bound.

\bigskip

\bigskip

We conclude with an example involving the realization of a spectrum with three
nonzero entries.

\bigskip

\begin{example}
 $\sigma=(\rho,\exp(\frac{\pi i}{10}),\exp(\frac{-\pi
i}{10}))$ has all its power sums positive if $\rho>\sqrt[9]{2\cos
(\pi/10)\text{ }}$ $=1.07....$. If we take $\rho=1.1$, and carry out the algorithm, we find that
$\sigma$ with $125$ zeros added is the spectrum of an $128\times128$ nonnegative matrix of the form of $X$ above. The least number of zeros required to be added to $\sigma$ to ensure realizability does not appear
to be known in this case.
\end{example}

\bigskip\bigskip

5. \ \ \ \textbf{References:}

\medskip

\bigskip\ 

[1] B.Beauzamy Products of polynomials and a priori estimates for coefficients
in polynomial decompositions.\ J. Symbolic Comput. 13 (1992) 463-472.

[2] A. Berman and R.J. Plemmons Nonnegative matrices in the Mathematical
Sciences. Second Edition \ SIAM 1994

[3] R. Bhatia, G. Krause and L. Elsner Bounds for the roots of a polynomial
and the eigenvalues of a matrix. Linear Algebra Appl. 142 (1990) 195-209.

[4] M. Boyle and D. Handelman The spectra of nonnegative matrices via symbolic dynamics.
 Ann. Math. 133 \ (1991) 249-316.\ \ 

[5] S. Friedland On an inverse problem for nonnegative and eventually
nonnegative matrices.
 Israel J. Math. 29 (1978) 43-60.

[6] S. Friedland A note on the nonzero spectrum of irreducible matrices.
arXiv: 0910.3415.

[7] A.Galantai and C.J. Hegedus Perturbation bounds for polynomials.
Numerische Math. 109 (2008) 77-100.

[8] C.R. Johnson Row stochastic matrices similar to doubly stochastic
matrices. Linear Multilinear Algebra 10 (1981) 113-130.

[9] T.J.Laffey and H. \v{S}migoc Nonnegative realization of spectra having
negative real parts. Linear Algebra\ Appl. 416 (2006) 148-159.

[10] T.J. Laffey and H. \v{S}migoc On a classic example in the nonnegative
inverse eigenvalue problem. ELA Electronic Journal of Linear Algebra 17 (2008) 333-342.

[11] R. Loewy and D. London \ A note on an inverse eigenvalue
problem for nonnegative matrices.\ Linear Multilinear Algebra 6 (1978) 83-90.

[12] M.E. Meehan Some results on matrix spectra. PhD thesis. National
University of Ireland, Dublin 1998.

[13] H.R. Suleimanova Stochastic matrices with real characteristic numbers. Dokl.Akad.NAUK,
SSSR (N.S.) 66 (1949) 343-345.

[14] \ J. Torre-Mayo, M.R. Abril-Raymundo, E. Alarcia-Estevez, C. Marijuan and
M. Pisanero The nonnegative inverse eigenvalue problem from the
coefficients of the characteristic polynomial :EBL Digraphs. Linear Algebra Appl.426 (2007) 729-773.

\ \ \ \ \ \ \ \ \ \ 

\ 
\end{document}